# The Application of Imperialist Competitive Algorithm for Fuzzy Random Portfolio Selection Problem

Mir Ehsan Hesam Sadati
Urmia University
Urmia, Iran

Jamshid Bagherzadeh Mohasefi
Urmia University
Urmia, Iran

## ABSTRACT
This paper presents an implementation of the Imperialist Competitive Algorithm (ICA) for solving the fuzzy random portfolio selection problem where the asset returns are represented by fuzzy random variables. Portfolio Optimization is an important research field in modern finance. By using the necessity-based model, fuzzy random variables reformulate to the linear programming and ICA will be designed to find the optimum solution. To show the efficiency of the proposed method, a numerical example illustrates the whole idea on implementation of ICA for fuzzy random portfolio selection problem.

## 1. INTRODUCTION
In many industries, there are many decision problems; i.e., scheduling problem, logistics. In these problems, it is important to predict future total returns and to decide an optimal asset allocation maximizing total profits under some constraints. It is easy to decide the most suitable allocation if we know future returns a prior. For earning the greatest profit, consideration of how to reduce a risk is utilized. We call such industrial assets allocation problems portfolio selection problems [1]. Markowitz [2] proposed the well-known mean–variance models and these models have been improved in both theory and algorithm [3, 4]. When selecting portfolio, an investor may encounter with both fuzziness and randomness. Fuzzy random variable can be a new useful approach to solve this kind of problem. A Fuzzy random variable was first introduced by Kwakernaak [5], and its mathematical basis was constructed by Puri and Ralescu [6]. In this paper, the asset return in portfolio selection problem are fuzzy random variables and we use the concept of necessity-based model to develop a solution method for the fuzzy random portfolio optimization problem.

This paper proposes the use of a new evolutionary algorithm known as Imperialist Competitive Algorithm (ICA) to solve fuzzy random portfolio optimization. The ICA is a meta-heuristic optimization method that is based on modeling for the attempts of countries to dominate other countries [7].

The rest of the paper is organized as follow: Section 2 includes basic concepts on fuzzy and fuzzy random theory. In section 3, the problem formulation is presented. Section 4, explains the ICA in detail. In section 5, a numerical example is solved to illustrate the proposed model. Finally, conclusion and future work will be presented in section 6.

## 2. BASIC CONCEPTS
In order to extend statistical analysis to situations when the outcomes of some random experiment are fuzzy sets, the concept of fuzzy random variable was introduced as an analogous notion to random variable. A Fuzzy random variable was first introduced by Kwakernaak [5], who introduced fuzzy random variables (FRV) as "random variables whose values are not real, but fuzzy numbers," and its mathematical basis was constructed by Puri and Ralescu [6]. Fuzzy random variables can be defined in an n dimensional Euclidean space $R^n$. We present the definition of a fuzzy random variable in a single dimensional Euclidean space R.

*Definition 1 [8]*

Let $(\Omega, A, P)$ be a probability space, where $\Omega$ is a sample space, A is a σ-field and P is a probability measure. Let $F_N$ be the set of all fuzzy numbers and B a Borel σ-field of R. Then a map $\tilde{\tilde{Z}} : \Omega \to F$ is called a fuzzy random variable if it holds that:

$$\left\{(\omega,\tau) \in \Omega \times R \mid \tau \in \tilde{\tilde{Z}}_\alpha(\omega)\right\} \in A \times B, \forall \alpha \in [0,1]$$

where

$$\tilde{\tilde{Z}}_\alpha(\omega) = \left[\tilde{\tilde{Z}}_\alpha^-(\omega), \tilde{\tilde{Z}}_\alpha^+(\omega)\right] = \left\{\tau \in R \mid \mu_{\tilde{Z}(\omega)}(\tau) \geq \alpha\right\}$$

is an α-level set of the fuzzy number $\tilde{\tilde{Z}}(\omega)$ for $\omega \in \Omega$.

*Definition 2*

LR fuzzy number $\tilde{A}$ is defined by following membership function:

$$\tilde{A}(X) = \begin{cases} L\left(\dfrac{A^0 - x}{\beta}\right) & if \quad A^0 - \beta \leq x \leq A^0 \\ 1 & if \quad A^0 \leq x \leq A^1 \\ R\left(\dfrac{x - A^1}{\gamma}\right) & if \quad A^1 \leq x \leq A^1 + \gamma \end{cases}$$

where $\left[A^0, A^1\right]$ show the peak of fuzzy number $\tilde{A}$ and $\beta, \gamma$ represent the left and right spread respectively; $L, R = [0,1] \to [0,1]$ with $L(0) = L(0) = 1$ and $L(1) = L(1) = 0$ are strictly decreasing, continuous functions. A possible representation of a LR fuzzy number is $\tilde{A} = \left(A^0, A^1, \beta, \gamma\right)_{LR}$.

## 3. FUZZY RANDOM PORTFOLIO SELECTION PROBLEM

In this paper, we deal with the following portfolio selection problem involving fuzzy random variable returns to maximize total future returns. In the following problem, called Fuzzy Random optimal portfolio selection problem, the return rate of assets are fuzzy random variables:

Problem 1

$$Max \ \tilde{\tilde{Z}} = \sum_{j=1}^{n} \tilde{\tilde{R}}_j x_j \quad (1)$$

$$s.t. \quad \sum_{j=1}^{n} x_j = M_0, \quad (2)$$

$$\sum_{j=1}^{n} \tilde{\tilde{R}}_j x_j \geq \tilde{\tilde{R}}_0 \quad (3)$$

$$0 \leq x_j \leq U_j \ ; \quad j = 1, 2, \ldots, n. \quad (4)$$

The parameters and variables are define as follow, for $j=1,2,\ldots,n$:

* $\tilde{\tilde{R}}_j = \left(\overline{R}_j^0, \overline{R}_j^1, \beta_j, \gamma_j\right)_{LR}$ represent Fuzzy random variables whose observed value for each $\omega \in \Omega$ is fuzzy number $\tilde{R}_j(\omega) = \left(R_j^0(\omega), R_j^1(\omega), \beta_j, \gamma_j\right)_{LR}$

* $\left(\overline{R}_j^0, \overline{R}_j^1\right) = \left(R_j^0 + \bar{t}R_j^2, R_j^1 + \bar{t}R_j^2\right)$ is a random vector which $\bar{t}$ is a random variable with normal distribution function $T$.

* $n$ = the number of assets for possible investment

* $M_0$ = available total fund

* $\tilde{\tilde{R}}_j$ = the rate of return of assets $j$ ( per period)

* $\tilde{\tilde{R}}_0$ = the return in dollars

* $x_j$ = decision variables which represent the dollar amount of fund invested in asset $j$

* $U_j$ = the upper bound of investment in asset $j$.

By using the necessity based-model, we reformulate the problem as a linear programming problem. In the next step, necessity-based model will be explained.

### 3.1 Necessity-Based Model

Since the obtain solution will be too optimistic, so necessity-based model can be suitable for pessimistic decision maker who wish to avoid risk. By Zadeh's extension [9] principle for objective function in problem 1, its membership function is given as follows for each $\omega \in \Omega$:

$$\mu_{\tilde{Z}(\omega)}(t) = \begin{cases} L\left(\dfrac{Z^0(\omega) - t}{\beta}\right) & if \quad t \leq Z^0(\omega) \\ 1 & if \quad Z^0(\omega) \leq t \leq Z^1(\omega) \\ R\left(\dfrac{t - Z^1(\omega)}{\gamma}\right) & if \quad otherwise \end{cases}$$

where

$$\tilde{Z}(\omega) = \left(Z^0(\omega), Z^1(\omega), \beta, \gamma\right)_{LR}, Z^0(\omega) = \sum_{j=1}^{n} R_j^0(\omega) x_j$$

, and $Z^1(\omega) = \sum_{j=1}^{n} R_j^1(\omega) x_j$.

The degree of necessity $N\left(\tilde{Z}(\omega) \geq f\right)$ under the possibility distribution $\mu_{\tilde{Z}(\omega)}(t)$ is defined as follows:

$$N\left(\tilde{Z}(\omega) \geq f\right) =$$
$$\inf_{y_1, y_2}\left\{\max\left\{1 - \mu_{\tilde{Z}(\omega)}(y_1), 1 - \mu_f(y_2)\right\} \mid y_1 \geq y_2\right\} \geq \eta$$

and the necessity degree of fuzzy constraint $\left(\sum_{j=1}^{n} \tilde{R}_j(\omega) x_j \geq \tilde{R}_0(\omega)\right)$ under the possibility distributions is defined as follows:

$$N\left(\sum_{j=1}^{n} \tilde{R}_j(\omega) x_j \geq \tilde{R}_0(\omega)\right) =$$
$$\inf_{y_1, y_2}\left\{\max\left\{1 - \mu_{\sum_{j=1}^{n} \tilde{R}_j(\omega) x_j}(y_1), 1 - \mu_{\tilde{R}_0(\omega)}(y_2)\right\} \mid y_1 \geq y_2\right\}$$

We maximize the degree of necessity $N\left(\tilde{Z}(\omega) \geq f\right)$ and the degree of necessity $N\left(\sum_{j=1}^{n} \tilde{R}_j(\omega) x_j \geq \tilde{R}_0(\omega)\right)$, our portfolio selection model in Problem 1 comes by the following model:

Problem 2

$$Max \ f \quad (5)$$

$$s.t. \quad \Pr\left\{\omega \mid N\left(\tilde{Z}(\omega) \geq f\right) \geq \eta\right\} \geq \lambda, \quad (6)$$

$$\sum_{j=1}^{n} x_j = M_0, \quad (7)$$

$$\Pr\left\{\omega \mid N\left(\sum_{j=1}^{n} \tilde{R}_j(\omega) x_j \geq \tilde{R}_0(\omega)\right) \geq \eta\right\} \geq \lambda, \quad (8)$$

$$0 \leq x_j \leq U_j \ ; \quad j = 1, 2, \ldots, n. \quad (9)$$

where $\lambda$ is a predetermined probability level and $\eta$ is a predetermined possibility level.

A feasible solution of portfolio selection problem is called a necessity solution. In order to transform the above model to a linear programming model, we need to reformulate (19) and (21). We use the Katagiri et al. [10] approach to linearization of above model. Consider the following theorem: [1]

*Theorem 1: [1, 10]*

$$\begin{cases} 1) \ \Pr\{\omega \mid N(\tilde{Z}(\omega) \geq f) \geq \eta\} \geq \lambda \Leftrightarrow \\ \sum_{j=1}^{n}\left(R_j^0 + T^*(1-\lambda)R_j^2\right)x_j - L^*(1-\eta)\sum_{j=1}^{n}\beta_j x_j \geq f \\ 2) \ \Pr\left\{\omega \mid N\left(\sum_{j=1}^{n}\tilde{R}_j(\omega)x_j \geq \tilde{R}_0(\omega)\right) \geq \eta\right\} \geq \lambda \Leftrightarrow \\ \sum_{j=1}^{n}\left(R_j^0 + T^*(1-\lambda)R_j^2\right)x_j - L^*(1-\eta)\sum_{j=1}^{n}\beta_j x_j \geq \\ R_0^0 + T^*(1-\lambda)R_0^2 - \beta_0 L^*(1-\eta) \end{cases}$$

where $T^*$, $L^*$ and $R^*$ are pseudo inverse functions defined as:

$$T^*(\lambda) = \inf\{t \mid T(t) \geq \lambda\}$$
$$L^*(\lambda) = \sup\{t \mid L(t) \geq \lambda\}$$
$$R^*(\lambda) = \sup\{t \mid R(t) \geq \lambda\}.$$

Now the optimal solution of Problem 2 is equal to the following linear parametric programming problem:

Problem 3

$$Max \quad \sum_{j=1}^{n}\left(R_j^0 + T^*(1-\lambda)R_j^2\right)x_j - L^*(1-\eta)\sum_{j=1}^{n}\beta_j x_j \tag{10}$$

$$s.t. \quad \sum_{j=1}^{n}x_j = M_0, \tag{11}$$

$$\sum_{j=1}^{n}\left(R_j^0 + T^*(1-\lambda)R_j^2\right)x_j - L^*(1-\eta)\sum_{j=1}^{n}\beta_j x_j \geq$$
$$R_0^0 + T^*(1-\lambda)R_0^2 - \beta_0 L^*(1-\eta) \tag{12}$$

$$0 \leq x_j \leq U_j; \quad j = 1, 2, ..., n. \tag{13}$$

Different methods have been proposed for solving portfolio optimization problem. Some of these methods are the computer simulation of such as genetic algorithms. In this paper, Imperialist Competitive Algorithm for portfolio optimization problem is proposed which is inspired by imperialistic competition. All the countries are divided into two types: imperialist states and colonies [11]. The following section will introduce the whole aspects of Imperialist Competitive Algorithm (ICA).

## 4. IMPERIALIST COMPETITIVE ALGORITHM

Imperialist competitive algorithm (ICA) is a new evolutionary optimization method that is inspired by imperialistic competition. It starts with an initial population like other evolutionary algorithms which is called country, which be of, type being colonized or being imperialist [11]. Imperialistic competition is the main part of proposed algorithm and causes the colonies to converge to the global minimum of cost function. This method is a new socio-politically motivated global search strategy that has recently been introduced for dealing with different optimization tasks [7].

Figure 1 shows the flow chart of the ICA [11].

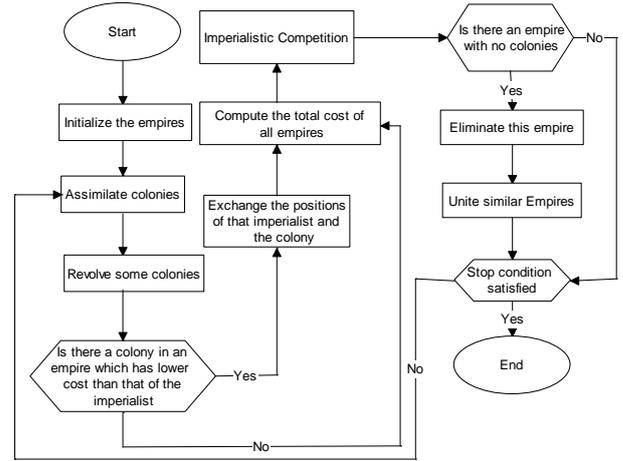

**Fig 1: The flow chart of ICA**

Some of the best countries in the population are selected to be the imperialist and the rest form the colonies of these imperialists. Each country is shown by $country_i = [p_1, p_2, ..., p_{N_{var}}]$ which $N_{var}$ is dimension of optimization problem [7]. The best countries in the initial population, considering the cost function of them, are selected as the imperialists and other counties known as the colonies of these imperialists.

The variable values in the country are represented as floating point number. The cost of a country is found by evaluating the cost function $f$ at the variables $(p_1, p_2, ..., p_{N_{var}})$. Therefore, we have [11]:

$$\cos t = f(country) = f(p_1, p_2, ..., p_{N_{var}}) \tag{14}$$

We select $N_{imp}$ of the most powerful countries to form the empires. The remaining $N_{col}$ of the population will be the colonies each of which belongs to an empire. To form the initial empires, we divide the colonies among imperialists based on their power. To divide the colonies among imperialists proportionally, the following normalized cost of an imperialist is defined:

$$C_n = c_n - \max(c_i) \tag{15}$$

Where $c_n$ is the cost of nth imperialist and $C_n$ is the normalized cost. As a result the normalized power of each imperialist can be defined:

$$p_n = \left|\frac{C_n}{\sum_{i=1}^{N_{imp}}C_i}\right| \tag{16}$$

After dividing colonies among imperialists, these colonies start closing to their empire. The total power of an empire can be determined by the power of imperialist country plus percentage of power of its colonies as follow [7]:

$$power_n = cost(imperialists_n) + \varepsilon \, mean\{(cost(colonies \, of \, empires_n))\} \quad (17)$$

Where $power_n$ is the total power (cost) of the *n*th empire and $\varepsilon$ is a positive number less than 0.1.

Powerless empires will collapse in the imperialistic competition and their colonies will be divided among other empires. After a while, all the empires except the most powerful one will collapse and all the colonies will be under the control of this empire [11], which is the answer of optimization problem.

## 4.1 Application of ICA to Portfolio Optimization Problem

For applying this algorithm to the portfolio optimization problem, Constraints handling should be done. In the initialization process, several initial solutions are created randomly for the ICA. To create a set of solutions, equality and inequality constraints should be satisfied. Handling the existing equality and inequality constraints in evolutionary optimization algorithm is very important. In this paper, we will use penalty function method [12] as follow to handle constraints in ICA.

### 4.1.1 Penalty Function Method [12]

Penalty functions are the oldest approach used to incorporate constraints into unconstrained optimization algorithms. The idea of this method is to transform a constrained optimization problem into an unconstrained one by adding a certain value to the objective function based on the amount of constraint violation present in a certain solution.

Consider our portfolio problem as:

Find $\vec{x}$ which optimizes $f(\vec{x})$

subject to:

$$g_i(\vec{x}) \leq 0, \quad i = 1,...,n$$
$$h_j(\vec{x}) = 0, \quad j = 1,...,p$$

The general formulation of exterior penalty function is:[12]

$$\phi(\vec{x}) = f(\vec{x}) \pm \left[\sum_{i=1}^{n} r_i \times G_i + \sum_{j=1}^{p} c_j \times L_j\right]$$

where $\phi(\vec{x})$ is the new objective function to be optimized, $G_i$ and $L_j$ are functions of the constraints $g_i(\vec{x})$ and $h_j(\vec{x})$, respectively, $r_i$ and $c_j$ are positive constants normally called "penalty factors".

The most common form of $G_i$ and $L_j$ is:

$$G_i = \max[0, g_i(\vec{x})]^\beta$$
$$L_j = |h_j(\vec{x})|^\gamma$$

where $\beta$ and $\gamma$ are normally 1 or 2.

Therefore, in our portfolio problem we can use this method for handling the both constraints (11, 12) to find the optimum solution by ICA.

## 5. AN EXAMPLE

In this section, an example is given to illustrate the proposed ICA through necessity-based model for portfolio optimization selection. Let us consider 5 securities whose returns are fuzzy random variables and their values are given in Table1.

The following parameters are used in ICA:
    Number of imperialist =10
    Number of countries = 100
    Revolution rate =0.2
    Maximum iteration =25
The parameters of other optimization methods are listed in [11].

**Table 1. Parameters of the example**

|   |   | $x_j$ |   |   |   |   |
|---|---|---|---|---|---|---|
| $j$ | 0 | 1 | 2 | 3 | 4 | 5 |
| $R_j^0$ | 250 | 1.3 | 1.2 | 1.35 | 1.4 | 1.45 |
| $R_j^1$ | 250 | 1.45 | 1.25 | 1.4 | 1.5 | 1.6 |
| $R_j^2$ | 50 | 0.6 | 0.5 | 0.5 | 0.6 | 0.6 |
| $\beta_j$ | 40 | 0.2 | 0.15 | 0.15 | 0.25 | 0.25 |
| $\gamma_j$ | 40 | 0.2 | 0.15 | 0.15 | 0.25 | 0.25 |

$\overline{t}$ is a normal random variable whose mean 0 and variance 1. The upper bound of investment amount in each stock is set to no more than 60 units of the total available fund. Given a total allocation budget of 200 units and annual return which is fuzzy random variable shown as $\tilde{\overline{R}} = M_0 \tilde{\overline{r}}_0$ where $\tilde{\overline{r}}_0 = (1.25 + 0.25\overline{t}, 1.25 + 0.25\overline{t}, 0.2, 0.2)$. Now we want to know what is the optimal solution for our portfolio selection problem for the different levels of probability and possibility {0.1, 0.4, 0.7, 0.9}. We apply the ICA to find the optimum solution of fuzzy random portfolio selection problem through the necessity-based model.

The best estimate solutions are obtained after 25 iterations using the Matlab program and optimum solutions according to the implementation of ICA are collected in Table 2.

**Table 2. Numerical results for Necessity-Based Model**

| $\lambda, \eta$ | $x_j$ | | | |
|---|---|---|---|---|
| | 0.1 | 0.4 | 0.7 | 0.99 |
| $x_1^*$ | 60 | 20 | 20 | 0 |
| $x_2^*$ | 0 | 0 | 0 | 60 |
| $x_3^*$ | 20 | 60 | 60 | 60 |
| $x_4^*$ | 60 | 60 | 60 | 20 |
| $x_5^*$ | 60 | 60 | 60 | 60 |
| OFV[a] | 422.54 | 289.3 | 187.48 | 95.56 |

[a] Objective function value.

According to the results and by Implementation of ICA for this problem, we obtained the best solutions and these results show the effectiveness of ICA for finding the optimum solution for portfolio optimization problem.

## 6. CONCLUSION

This paper proposed a solution method for portfolio selection model whose parameters were fuzzy random variables. The idea was based on necessity-based model and the way of finding the optimum solution was imperialist competitive algorithm. Each individual of the population was called a country. The population was divided into two groups: colonies and imperialist states. Penalty function method for handling the equality and inequality constraints was utilized which always provide feasible solution. For future research, we will apply the other methods and algorithms to find the optimum solution for fuzzy random portfolio selection.

## 7. REFERENCES

[1] M. E. Hesam sadati, J. Nematian, " Two-level linear programming for fuzzy random portfolio optimization through possibility and necessity-based model", Procedia Economics and finance 5 (2013) 657-666

[2] H. Markowitz, "Portfolio selection", Journal of Finance 7 (1952) 77–91.

[3] Y. Crama, M. Schyns, "Simulated annealing for complex portfolio selectitn problems", European Journal of Operational Research 150 (2003) 546–571.

[4] Y. Xia, B. Liu, S. Wang, K. lai, "A model for portfolio selection with order of expected returns", Computers and Operations Research 27 (2000) 409–422.

[5] H. Kwakernaak, "Fuzzy random variables", Definitions and theorems, Inf. Sci 15 (1978) 1.

[6] M.L. Puri, D.A. Ralescu, "Fuzzy random variables", Journal of Math and Anal 114 (1986) 409.

[7] G. Mokhtari, et al. "Application of Imperialist Competitive Algorithm to Solve Constrained Economic Dispatch", International Journal on Elctrical Engineering and Informatics, 4 (2012) 4

[8] M. Sakawa, "Fuzzy Sets and Interactive Multi-objective Optimization", Plenum Press, New York (1993)

[9] L. A. Zadeh, "Fuzzy sets.Informations and controls", 8 (1965) 338–353.

[10] H. Katagiri, M. Sakawa, K. Kato, "Interactive multi-objective fuzzy random linear programming: maximization of possibility and probability", European journal of operational research 188 (2008) 530.

[11] E. Atashpaz-Gargari, C. Lucas, "Imperialist Competitive Algorithm: An Algorithm for Optimization Inspired by Imperialistic Competition", IEEE Congress on Evolutionary Computation, (2007)

[12] C. A. Coello Coello, " Theoretical and Numerical Constraint-Handling Techniques used with Evolutionary Algorithms: A Survey of the State of the Art", Computer Methods in Applied Mechanics and Engineering, 191 (2002) 1245-1278